\input amstex
\documentstyle{amsppt}
\NoBlackBoxes
\magnification 1200
\hsize 6.5truein
\vsize 8truein
\topmatter
\title Kakeya sets in Cantor directions \endtitle
\author  Michael D. Bateman and Nets Hawk Katz \endauthor
\affil Indiana University \endaffil
 \thanks The first author was supported by NSF grant DMS-FRG-0139874;
  The second  author was supported by NSF grant DMS 0432237.
\endthanks 
\endtopmatter

\head \S 0 Introduction \endhead

In this paper, we prove the following.

\proclaim{Theorem 0.1} For any $N=3^n$, there is a union of $N$ parallelograms $P_1, \dots P_N$
in $\Bbb R^2$ of eccentricity $\sim N$ and area $\sim {1 \over N}$ so that the slopes of the
long sides of $P_1,\dots,P_N$ are all contained in the standard middle-thirds Cantor set, so
that 
$$| \bigcup_{j=1}^N P_j | \lesssim {1 \over \log N}, \tag 0.1 $$
but so that, if we let $2P_j$ be the double of the parallelogram we have
$$| \bigcup_{j=1}^N 2P_j | \gtrsim {\log \log N \over \log N}.$$
\endproclaim

In the statement of the theorem as in the rest of the paper, we use the convention that
when $S$ is a subset of the plane $\Bbb R^2$, we denote by $|S|$ the Lebesgue measure
of $S$. Further when $A_N$ and $B_N$ are numbers depending on $N$, and we write
$A_N \lesssim B_N$, we mean there is a constant $C$ independent of $N$ so that
$$A_N \leq C B_N.$$

The proof of our theorem is by a probabilistic construction. The estimate which allows
us to prove inequality (0.1) is a  now fairly standard estimate on percolation on trees following
the work of Russ Lyons ([L],[LP]). As far as we know, this idea has not appeared in the
study of Kakeya sets before. The moral of the story is that if we define (loosely) a Kakeya set in the
plane as a ``1 dimensional" family of unit line segments whose union has measure 0 then
while it is true that the random family of line segments is not a Kakeya set, it is the case that
the random, {\it sticky}, set of line segments is a Kakeya set. Here we use the term sticky as in [KLT].

If we let $\Cal S$ be the set of line segments in the plane whose slope is in the standard 
Cantor set, and we define for $s \in \Cal S$, the expression $av_s(f)$ to be the average of
a function $f$ on $s$, where $f$ must be locally integrable on lines, we may define
a maximal operator
$$\Cal M f(x) = \sup_{x \in s \in \Cal S} av_s |f|.$$
An immediate consequence of our theorem is

\proclaim{Corollary 0.2} The maximal operator $\Cal M$ is unbounded on any $L^p(\Bbb R^2)$
with $p \neq \infty$. \endproclaim

This was proved for $p \leq 2$ in [K]. Previously, the operator had been explicitly
studied in [DV] and [V]. The boundedness of this operator had been known in the
folklore as an open problem for more than a decade previously.

The proof in [K] had no applications to Fourier multiplier operators
precisely because Theorem 0.1 had not been proved.  We intend to address the
implications of our present result for multiplier operators and for the theory of
directional maximal functions, in general, in future work.

We thank Russ Lyons for helpful discussions.

\head \S 1  Geometric Constructions \endhead

We denote by $T_n$ the set of all $n$-digit strings $.a_1a_2 \dots a_n$ with each
$a_j$ taking on the value $0,1$ or $2$. Here we consider $T_0$ to be the
singleton set containing ``$.$", the empty decimal.  We define the maps
$$\pi_j : T_n \longrightarrow \{0,1,2\},$$
by
$$\pi_j(.a_1a_2 \dots a_n)=a_j,$$
and for $j < n$, we define
$$\pi^j: T_n \longrightarrow T_j,$$
by
$$\pi^j(.a_1a_2 \dots a_n)=.a_1a_2 \dots a_j.$$
We define
$$T_n^*=\bigcup_{j=0}^n T_j.$$
We may view $T_n^*$ as a rooted ternary tree with an edge between
$s \in T_j$ and $t \in T_{j-1}$ whenever $\pi^{j-1}(s)=t$.
We denote this edge by $e_{t,\pi_j(s)}$ and say that $s$ is
the $\pi_j(s)$th child of $t$. We identify the tree $T_n^*$ with the
triadic intervals of length greater than $3^{-n}$, by the map
$$I(s)=[s,s + {1 \over 3^j}],$$
when $s \in T_j$ and $s=.a_1 \dots a_j$ is identified with the triadic rational
$${a_1 \over 3} + {a_2 \over 9} + \dots + {a_j \over 3^j}.$$
Whenever $s,t \in T_n^*$, and $I(s) \subset I(t)$, we say that
$t$ is an ancestor of $s$, or $s$ is a descendant of $t$.

We denote by $C_n \subset T_n$ the set of all $n$-digit strings $.a_1a_2 \dots a_n$
so that each $a_j$ takes on the value either $0$ or $2$.  Then $I(C_n)$ is the $n$th
stage of the construction of the standard Cantor set. We will say that a map
$$\sigma: T_n \longrightarrow C_n,$$
is {\it sticky} provided that for any $s \in T_n$, the value of
$\pi_j (\sigma(s))$ depends only on $\pi^j(s)$. We shall define a random variable
$\sigma_n$ which takes values in sticky maps from $T_n$ to $C_n$. This random
variable shall in fact be evenly distributed among such maps, but we define its components
more explicitly.

To each edge $e_{t,a}$ of $T_n^*$, we define a random variable $r_{t,a}$. The variables
$r_{t,a}$ are independent and take on the values $0$ and $2$ with probability ${1 \over 2}$
each. Now we define
$$\sigma_n(s)=s^{\prime},$$
where $\pi_j(s^{\prime})=r_{\pi^{j-1}(s),\pi_j(s)}$.

Following [K], we assign a ``Kakeya set" to every possible value of the random variable
$\sigma_n$. (In [K], this was actually when $r_{t,a}=0$ for $a=0,1$ and $r(t,a)=2$ for $a=2$,
independently of $t$.)  Given a sticky map
$$\sigma:T_n \longrightarrow C_n,$$ 
we define for each $s \in T_n$, a parallelogram in $\Bbb R^2$ which we will denote by
$P_{\sigma,s}$. The parallelogram $P_{\sigma,s}$ has as its corners the points
$(0,{s \over 3}),(0,{s \over 3}+{{1 \over 3^{n+1}}}), (1, {s \over 3} + \sigma(s)),$
and $(1,{s \over 3}+{{1 \over 3^{n+1}}} + \sigma(s))$. (Here we again identify
$s$ and $\sigma(s)$ as real numbers by the ternary expansion.) We think
of $P_{\sigma,s}$ as a tube with eccentricity approximately ${1 \over 3^{n+1}}$ which
begins at $(0,{s \over 3})$ and has slope $\sigma(s)$. Then we define a ``Kakeya set"
by
$$K_{\sigma} =\bigcup_{s \in T_n} P_{\sigma,s}.$$

Our first goal is to prove

\proclaim{Lemma 1.1} For any choice of a sticky map
$$\sigma: T_n \longrightarrow C_n,$$
we have that 
$$|K_{\sigma}| \gtrsim {\log n \over n}.$$
\endproclaim

Notice that Lemma 1.1 is a generalization of ([K],Lemma 2.3).

To prove this, we first establish the following elementary uniformity inequality in measure theory.

\proclaim{Proposition 1.2} Suppose $(X,\Cal M,\mu)$ is a  measure space and
$A_1,\dots,A_K$ are sets with $\mu(A_j)=\alpha$. Let $m>0$. Suppose that
$$\sum_{j=1}^K \sum_{k=1}^K \mu(A_j \cap A_k) \leq K m \alpha,$$
then
$$\mu(\bigcup_{j=1}^K A_j) \geq {K \alpha \over 16m}.$$
\endproclaim

(The 16 in the denominator is unnecessary, but simplifies the proof slightly.)

\demo{Proof} It must be that there is $S \subset \{1, \dots K\}$ with
$\#(S)   \geq{K \over 2}$ so that we have
$$\sum_{j=1} \mu(A_j \cap A_k) \leq  2m \alpha,$$
whenever $k \in S$. For any such $k$, there must be a measurable set 
$B_k \subset A_k$ so that
$$\sum_{j=1}^K \chi_{A_j} (x) \leq 4m,$$
for any $x \in B_k$, and so that $\mu(B_k) \geq {\alpha \over 2}$.
Then
$$\int \sum_{k \in S} \chi_{B_k} \geq {K \alpha \over 4},$$
but 
$$\sum_{k \in S} \chi_{B_k}(x) \leq \sum_{j=1}^K \chi_{A_j}(x) \leq 4m,$$
for $x \in \cup_{k \in S}  B_k$.
Thus by Chebychev's inequality, we have
$$ {K \alpha \over 16 m} \leq \mu( \bigcup_{k \in S} B_k) \leq \mu(\bigcup_{j=1}^K A_j),$$
which was to be shown. \qed \enddemo

\demo{Proof of Lemma 1.1} We will show that for $0 \leq j < \log n$, with
$S_j = [3^{-j},3^{1-j}] \times \Bbb R$
we have the estimate
$$|K_{\sigma} \cap S_j| \gtrsim {1 \over n}.$$
We see that
$$K_{\sigma} \cap S_j = \bigcup_{s \in T_n} P_{\sigma,s,j},$$
where
$$P_{\sigma,s,j} = P_{\sigma,s} \cap S_j.$$
Since for each value of $s$, we have
$$|P_{\sigma,s,j}| ={2 \over 3^{j+n}},$$
it suffices to show, in light of Proposition 1.2 that
$$\sum_{s_1 \in T_n} \sum_{s_2 \in T_n} 
|P_{\sigma,s_1,j} \cap P_{\sigma,s_2,j}|  \lesssim {n \over 3^{2j}}. \tag 1.1$$
(Note that the inequality fails for $j \geq \log n$  because of the diagonal
part of the sum.)
Between any $s_1$ and $s_2$ we define the triadic distance $d(s_1,s_2)$ to be $3^{-k}$
where $k$ is the largest number for which $\pi^k(s_1)=\pi^k(s_2)$. Note that for any $s_1 \neq s_2$,
we have that 
$$P_{\sigma,s_1,j} \cap P_{\sigma,s_2,j} \neq \emptyset$$
implies that $d(s_1,s_2) \gtrsim 3^j |s_1 -s_2|,$ where again we have identified $s_1$ and $s_2$ as
numbers. Further, we always have the estimate
$$|P_{\sigma,s_1,j} \cap P_{\sigma,s_2,j}| \lesssim  {1\over 3^{n+j} |s_1-s_2|} , \tag 1.2$$
because $3^j |s_1-s_2|$, bounds below the difference in the slopes $\sigma(s_1)$ and $\sigma(s_2)$.
We divide up the sum in (1.1) according to the approximate value of $3^j |s_1-s_2|$ and
observe that  letting $A_{k,j}$ be
the number of pairs $(s_1,s_2)$ for which
$d(s_1,s_2) \geq 3^j |s_1-s_2|$ and  $3^j |s_1-s_2| \sim 3^k$, we have
$$A_{k,j} \lesssim  3^{n+k-2j}. \tag 1.3$$
Combining (1.2) and (1.3) and summing over $k$ proves the estimate (1.1).
\qed \enddemo

For the remainder of this section, we fix a point $(t,y) \in \Bbb R^2$ with
${1 \over 3} < t < 1$. We investigate the probability $P_n(t,y)$ 
of the event that $(t,y) \in K_{\sigma_n}$.

For every $s \in T_k$ and every $c \in C_k$, we consider $I_{s,c,t}$
which is the set of $y$ so that $(t,y)$ is contained in a line whose $y$-intercept
is in the interval $[{s \over 3}, {s \over 3} + {1 \over 3^{k+1}}]$ and whose
slope is contained in $[c, c +{1 \over 3^k}]$. We easily see that
$$I_{s,c,t}= [{s \over 3} + tc, {s \over 3} + {1 + 3t \over 3^{k+1}}  + tc ].$$
We observe that for any distinct $c_1,c_2 \in C_k$, we have
$|c_1-c_2| \geq {2 \over 3^k},$ so that since $t > {1 \over 3}$,
the collection
$$\{ I_{s,c,t} \}_{c \in C_k},$$
is pairwise disjoint. Therefore for each value of $s$, there is at most
one value of $c$ so that $y \in I_{s,c,t}$. (There may be no such value.) If such
a value $c$ exists we denote it by $c=c_{t,y}(s)$. Otherwise, we write
$c_{t,y}(s)=\infty$. Note that, by definition, if $c_{t,y}(s)$ is finite then $c_{t,y}(s^{\prime})$ is
finite for any ancestor $s^{\prime}$ of $s$. Note further that if we are given $s_1$ and $s_2$ with
$c_{t,y}(s_1),c_{t,y}(s_2)$ both finite and if $I(s_2) \subset I(s_1)$ then
$I(c_{t,y}(s_2)) \subset I(c_{t,y}(s_1))$. 
We denote by $T^*_{n,t,y}$, the set of those
$s \in T^*_n$ so that $c_{t,y}(s)$ is finite. Then the collection $T^*_{n,t,y}$
is a subtree of $T^*_n$.

We make two observations about the tree $T^*_{n,t,y}$. The first observation is that
the event $(t,y) \in K_{\sigma_n}$ occurs only if there is some $s \in T_n \cap T^*_{n,t,y}$
so that $\sigma_n(s)=c_{t,y}(s)$. This, in turn, happens if and only if for every $0<k \leq n$
we have that
$$\pi_k(c_{t,y}(s)) =r_{\pi^{k-1}(s),\pi_k(s)}. \tag 1.4,$$
The events in (1.4) are in one to one correspondence with the edges
$e_{\pi^{k-1}(s),\pi_k(s)}$, are independent of one another, and occur with probability
${1 \over 2}$. Thus $P_n(t,y)$ is bounded by the probability that if we remove each edge
of $T^*_{n,t,y}$ independently with probability ${1 \over 2}$, that we leave in
place a path from the root to the $n$th generation. This is called, in the probability literature,
(see {\it e.g.} [L],[LP]) the survival probability of Bernoulli( ${1 \over 2}$) percolation
on the tree $T^*_{n,t,y}$. We record this observation as a Lemma.

\proclaim{Lemma 1.3} With ${1 \over 3} < t < 1$, we have that
$P_n(t,y)$, the probability that $(t,y)$ is in the random ``Kakeya set"
$K_{\sigma_n}$ is bounded by the survival probability of
Bernoulli(${1 \over 2}$) percolation on the associated tree
$T^*_{n,t,y}$.
\endproclaim

The second observation is that for any $k$, the set of $s \in T_n$ such that $y \in I_{s,c,t}$ is contained
in $2^k$ intervals of length $t 3^{-k}$ which in turn is contained in $\lesssim 2^k$
triadic intervals of length $3^{-(k+1)}$. Thus, we get immediately

\proclaim{Lemma 1.4} We have, for every $0 \leq k \leq n$ the estimate
$$\#(T_k \cap T^*_{n,t,y}) \lesssim 2^k.$$
\endproclaim

Lemmas 1.3 and 1.4 will be enough to allow us to obtain the estimate which
we require for $P_n(t,y)$. We carry this out in the following section.

\head \S 2 Percolation on Trees \endhead

In this section, we review part of the theory of percolation on trees. We do not claim
any originality. All results are special cases of theorems of Russ Lyons (see {\it e.g.} [L],[LP]
. Pointers may be found there to a much wider literature). We thank Russ Lyons for
explaining his work to us.

We let $T^{\prime} \subset T^*_n$ be a subtree. We remove each edge of $T^{\prime}$ independently
with probability ${1 \over 2}$. We denote by $P(T^{\prime})$ the probability that a path remains
from the root to $T_n \cap T^{\prime}$.

We introduce one other quantity associated to $T^{\prime}$. We view $T^{\prime}$ as an electric
circuit which has a battery whose positive node is connected to the root and whose
negative part is connected in parallel to each vertex of $T_n \cap T^{\prime}$. On each edge
of $T^{\prime}$ which connects a vertex of $T_{k-1}$ to a vertex of $T_k$, we place
a resistor with resistance $2^k$. We denote by $R(T^{\prime})$, the resistance between the root
of $T^{\prime}$ and the bottom $T^{\prime} \cap T_n$. (For more on the mathematical
theory of electrical circuits, see [LP].) The following theorem is due to Lyons [L],  in greater
generality and with a better constant. We include the proof which follows simply to make the paper
self-contained.

\proclaim{Theorem 2.1 (Lyons)} We have that
$$P(T^{\prime}) \lesssim {1 \over 2 + R(T^{\prime})}.$$
\endproclaim

\demo{Proof} We prove this by induction on $n$. Clearly it is true for constant 2, when $n=0$.
We assume up to $n-1$, we have
$$P(T^{\prime}) \leq {12 \over 2 + R(T^{\prime})}.$$

We observe that if $T^{\prime}$ is subtree of $T$ containing the root, we may view
$T^{\prime}$ as the root, together with up to 3 edges connected to 3 trees $T_1,T_2,$ and $T_3$.
(If some of these trees are empty, we assign them probabilty zero and infinite resistance.) 
We denote
$$P(T_j)=P_j,$$
and
$$R(T_j)=R_j.$$
Then
we have the recursive formulae
$$P(T)={1 \over 2} (P_1 + P_2 + P_3) - {1 \over 4} (P_1P_2 + P_1P_3 + P_2P_3)
+{1 \over 8} P_1 P_2 P_3 \tag 2.1$$
and
$${1 \over R(T)}={1 \over 2+2R_1} + {1 \over 2 + 2R_2}
+{1 \over 2 + 2 R_3}. \tag 2.2$$

Now we break into two cases. In the first case, we have ${12 \over 2 + R_j} > 2$ for
some $j$. Then we have $R_j < 4$. This implies $R(T)<10$ which implies 
${12 \over 2 + R(T)} > 1$, so that we certainly have 
$$P(T) \leq {12 \over 2 + R(T)}.$$

We define 
$$Q_j={12 \over 2 + R_j}.$$
We may assume each $Q_j \leq 2$.
Observe that if we define
$$F(x,y,z)= 1 -(1-{1 \over 2} x)(1-{1 \over 2}y)(1 -{1 \over 2}z),$$
on the domain $[0,2] \times [0,2] \times [0,2]$ then $F$ is monotone increasing in
each variable.
Therefore we have that
$$\eqalign{ P(T) &=  F(P_1,P_2,P_3) \cr
                              &\leq  F(Q_1,Q_2,Q_3) \cr
                              &\leq {1 \over 2}(Q_1+Q_2+Q_3) - {1 \over 6} (Q_1Q_2 + Q_1 Q_3 +Q_2 Q_3) }.\tag 2.3$$
Note that the equality is (2.1), while for the two inequalities we have used that the $Q$'s are 
$\leq 2$.

Now plugging into (2.3), the definition of the $Q$'s,  we obtain
$$\eqalign{ P(T) 
&\leq {12 \over 2} [ {(R_1+2)(R_2+2) + (R_1+2)(R_3+2) + (R_2+2)(R_3+2) - {12  \over 6} 
(R_1+R_2+R_3 + 6) \over (R_1+2)(R_2+2)(R_3+2) }] \cr
&\leq {12\over 2} [ {(R_1+2)(R_2+2) + (R_1+2)(R_3+2) + (R_2+2)(R_3+2) - {12 \over 6} 
(R_1+R_2+R_3 + 6) \over (R_1+2)(R_2+2)(R_3+2) - 4R_1-4R_2-4R_3-13}] \cr
&\leq {12 \over 2} [ {(R_1+1)(R_2+1) + (R_1+1)(R_3+1) + (R_2+1)(R_3+1)
\over (R_1+2)(R_2+2)(R_3+2) -4R_1-4R_2-4R_3-13}] \cr
&= {12 \over R(T) +2}. } $$
Here the second inequality is by
decreasing the denominator and the third inequality is by increasing the numerator.

\qed \enddemo

Next we estimate the resistance of the trees we are interested in.

\proclaim{Lemma 2.2} Let $T^*_{n,t,y}$ be as in section 1. Then
$$R(T^*_{n,t,y}) \gtrsim n.$$ \endproclaim

\demo{Proof} We use the basic physical principle, that the resistance of any circuit
may be reduced by shortcircuiting it with perfect conductors. We identify all vertices in each $T_k$,
thus reducing the resistance.
Then by Lemma 1.4, we have that $T_{k-1}$ and $T_k$ are connected by $\lesssim 2^k$ resistors
of resistance $2^k$ connected in parallel. Thus the resistance between $T_{k-1}$ and $T_k$
is $\gtrsim 1$. Thus the total resistance is $\gtrsim n$. \qed \enddemo

\proclaim{Corollary 2.3} Let ${1 \over 3} < t \leq 1$. Then with $P_n(t,y)$, the probability
that $(t,y) \in K_{\sigma_n}$, we have that
$$P_n(t,y) \lesssim {1 \over n} .$$\endproclaim

\demo{Proof} We combine Lemma 1.3, Theorem 2.1, and Lemma 2.2. \qed \enddemo

\head \S 3 Proof of the main theorem \endhead

\demo{Proof} We observe that in order for a point $(t,y)$ to be in any set  $K_{\sigma}$,
it must be that $0 \leq y \leq {4 \over 3}$. Thus $E$, the expected measure of
$K_{\sigma_n} \cap ([{1 \over 3},1] \times \Bbb R)$ is given by
$$E=\int   (\int_{1 \over 3}^1 \int_0^{4 \over 3} \chi_{K_{\sigma}} (t,y) dy dt) d\sigma,$$
where the outside integral takes place on a finite probability space. Interchanging the
integrals, we see that
$$E=\int_{1 \over 3}^1 \int_0^{4 \over 3} P_n(t,y) dy dt \lesssim {1 \over n}.$$
Therefore there is a choice of $\sigma$ for which 
$$|K_{\sigma} \cap ([{1 \over 3},1] \times \Bbb R)| \lesssim {1 \over n}. $$
On the other hand
$$|K_{\sigma}| \gtrsim {\log n \over n}.$$
Thus $K_{\sigma}$ is the desired example.

\end
\Refs\nofrills{References}

\widestnumber\key{KLT}


\ref \key DV \by Duoandikoetxea, J. and Vargas, A. \paper Directional Operators
and radial functions on the plane \jour Ark. Mat. \vol 33 \yr 1995 \pages 281-291 \endref

\ref \key L \by Lyons, R. \paper Random walks, Capacity, and Percolation on trees
\jour Ann. Probab. \vol 20 \yr 1992 \pages 2043-2088 \endref

\ref \key LP \by Lyons, R. and Peres, Y. \book Probability on Trees and Networks, in
preparation, http://mypage.iu.edu/~rdlyons/ $\#$ book \endref



\ref \key K \by Katz, N.H.
\paper A counterexample for maximal operators over a Cantor set of directions
\jour Mat Res. Let.
 \vol 3 \yr 1996 \pages  527--536 \endref
 
 \ref \key KLT \by Katz, N.H, Laba, I., and Tao, T.
 \paper An improved bound on the Minkowski dimension of Besicovitch sets
 \jour Ann. Math. \vol 152 \yr 2000 \pages 383--446 \endref
 
 \ref \key V \by Vargas, A. \paper A remark on a maximal function over a Cantor set of directions
 \jour Rend. Circ. Mat. Palermo \vol 44 \yr 1995 \pages 273--282 \endref








 
 \endRefs
